\documentclass{amsart}
\newtheorem{lemma}{Lemma}
\newtheorem{theorem}{Theorem}

\newtheorem{definition}{Definition}

\newcommand{\eps}{\varepsilon}
\newcommand{\R}{\mathbb{R}}

\newcommand{\overeq}[1]

\newcommand{\beqar}{\begin{eqnarray*}}
\newcommand{\eeqar}{\end{eqnarray*}}

\newcommand{\beqarl}{\begin{eqnarray}}
\newcommand{\eeqarl}{\end{eqnarray}}

\newcommand{\be}{\begin{equation}}
\newcommand{\ee}{\end{equation}}

\def\p{\partial}

\def\R{\mathbb{R}}

\usepackage[english]{babel} % English formatting
\usepackage[utf8]{inputenc} % Standard encoding
\usepackage[a4paper,left=2.9cm,right=2.9cm,top=3cm,bottom=3cm]{geometry} % Page formatting
\usepackage{indentfirst} % Indents the first paragraph
\usepackage{amsmath} % Maths type package
\usepackage[dvipsnames]{xcolor}
\usepackage{bm} % Bold font maths
\usepackage{graphicx} % Advanced graphics package
\usepackage{amssymb}

\usepackage{comment}
\usepackage{tikz-cd}
\usetikzlibrary{positioning, calc}
\usepackage{xcolor}
\usepackage[export]{adjustbox} 
\usepackage{pdflscape} % Make pages landscape
\usepackage{fancyhdr} % Fancy headers
\usepackage[colorlinks=true,citecolor=blue,urlcolor=blue,linkcolor=violet]{hyperref} % Link colours
\usepackage{amsfonts}
\usepackage{flafter} % Reference any 'float'
\usepackage[framemethod=tikz]{mdframed} % Box off stuff
\usepackage{xcolor}
\usepackage{mathtools}
\usepackage{color} % Colour support
\usepackage{wrapfig} % Text flowing around figures
\usepackage{lipsum} % Generates meaningless text
\usepackage{url}
\graphicspath{ {Images/} } % Folder for images 
\usepackage[makeroom]{cancel}
\usepackage{multicol} % Multi Column Environment
\usepackage{float} % Lets you add images to multi-column environment

\usepackage[numbers]{natbib}

\title{First order non-instantaneous corrections in collisional kinetic alignment models}

\author{Laura Kanzler}
\email{laura.kanzler@sorbonne-universite.fr}
\address{Sorbonne Universit\'e, UMR CNRS 7598, Universit\'e de Paris Cit\'e, Laboratoire Jacques-Louis Lions, 75005 Paris, France.}

\author{Carmela Moschella}
\email{carmela.moschella@univie.ac.at}
\address{University of Vienna, Faculty for Mathematics, Oskar-Morgenstern-Platz 1, 1090 Wien, Austria.}

\author{Christian Schmeiser}
\email{christian.schmeiser@univie.ac.at}
\address{University of Vienna, Faculty for Mathematics, Oskar-Morgenstern-Platz 1, 1090 Wien, Austria.}

\begin{document}

\pagenumbering{gobble} % keep title page without a number
\begin{abstract}
In this work the standard kinetic theory assumption of instantaneous collisions is lifted. A model for higher order non-instantaneous
alignment collisions is presented and studied in the asymptotic regime of short collision duration. A first order accurate approximative model is derived as a correction to the
instantaneous limit. Rigorous results on its well-posedness and on the instantaneous limit are proven. The approximative model is a system of two equations. An equally accurate 
scalar approximation is suggested.
\end{abstract}


\maketitle

\pagenumbering{arabic} % start page numbers again
\setcounter{page}{1}

\section{Introduction}
The Boltzmann equation of gas dynamics \cite{cercignani2013mathematical} is based on the simplifying assumption that collisions between particles are hard, i.e., instantaneous, 
such that the particle dynamics in phase space is governed by a velocity jump process. The same is true for various kinetic models for living agents like bacteria 
\cite{aranson2005pattern}, \cite{bertin2009hydrodynamic}, \cite{calvez2014confinement}, \cite{carlen2015boltzmann}, \cite{harvey2011study}, \cite{hittmeir2020kinetic}, 
\cite{MURPHY2024109266}, undergoing
spontaneous velocity changes or hard collisions, but also for models of opinion formation, with instantaneous changes of opinion \cite{during2015opinion}, \cite{toscani2006kinetic}. 
In gas dynamics, the hard collision assumption also justifies 
the restriction to binary collisions, since collisions of more than two particles are too rare to have an influence on the particle distribution \cite{cercignani2013mathematical}.
There have been efforts, however, to extend the Boltzmann equation to also include three-particle-collisions \cite{ampatzoglou2021rigorous}, \cite{sengers1973three}.
Non-instantaneous collisions have apparently only been addressed in the context of quantum particles  \cite{lipavsky1999noninstantaneous}.

In \cite{kanzler2022kinetic} the authors have started an investigation of kinetic models with non-instantaneous collisions considering a model problem for particle alignment, where 
in binary collision processes of positive duration, the one-dimensional velocity variables of pairs of particles approach their mean value. This can be seen as a version of the 
Vicsek model \cite{vicsek1995novel} (see \cite{bolley2012mean}, \cite{briant2022cauchy} for kinetic formulations), where the interaction is only pairwise and it is turned on and off
stochastically. Two versions of the model have been 
considered: one where collision processes have deterministic duration and end after the mean value (i.e., complete alignment) has been reached; and another one with stochastic 
collision duration governed by a Poisson process. The well posedness of the model problems has been studied in \cite{kanzler2022kinetic} as well as their instantaneous limits as 
the collision duration tends to zero. The instantaneous limit problems are hard collision kinetic models of standard form. The long time limit is a fully aligned state, where the 
distribution function collapses to a Delta distribution. This is a consequence of the energy loss in the collision processes, a property this model shares with other models for alignment 
\cite{aranson2005pattern}, \cite{carlen2015boltzmann}, \cite{hittmeir2020kinetic}, \cite{MURPHY2024109266}, and with the inelastic Boltzmann equation \cite{bobylev2000some}, \cite{carrillo2007contractive}, \cite{mischler2006cooling}.

The present work can be seen as a continuation of \cite{kanzler2022kinetic}. It starts from a model including higher order non-instantaneous collisions, where more than two particles 
interact. The model, presented in the following section, takes the form of a system of coagulation-fragmentation equations \cite{ball1990discrete} with additional drift terms. 
Coagulation and fragmentation correspond to (groups of) particles joining and, respectively, leaving a collision process, whose internal dynamics is described by the drift. The main 
goal is to consider the situation of short collision duration and to find (first order) corrections to the instantaneous limit model.

The following section has three parts. In the first part the higher order non-instantaneous collision model is presented. The second part is devoted to the discussion of formal properties, leading to a conjecture on the long time behavior.
This part is not essential for the further development of this work. The third part contains the main modeling step. It starts with a formal discussion of the instantaneous limit producing a kinetic model of standard form with binary instantaneous collisions. Then first order corrections are added resulting in a system of two equations for the distribution of free particles between collisions and for the distribution of pairs of particles involved in binary collision processes. These equations also contain an account of three-particle-collisions.
 
Rigorous results on this system are contained in Section 3. We prove an existence and uniqueness result for mild solutions as well as a rigorous justification of the instantaneous
 limit. Section 4 is devoted to the question of finding a formally equally accurate approximative model, which can be written as a scalar equation for a one-particle distribution.
 It has already been noted in \cite{kanzler2022kinetic} that the model for non-instantaneous binary collisions can be written as a scalar equation with time delays, which are small
 for small collision duration. It is then rather straightforward to derive first order corrections by Taylor expansion \cite{driver1973equation}, \cite{kurzweil2006small}. 
 Unfortunately, this asymptotic approximation is not structure preserving. For example, the resulting scalar equation does not preserve the nonnegativity of the solution in general.  
 Therefore we propose a scalar model with delays, which is both formally accurate up to first order and preserves nonnegativity of the solution.

We conclude the introduction with an overview of the various models presented in the following, in decreasing order of complexity:

\begin{center}
\begin{figure}[H]
\begin{tikzpicture}[scale=0.5,
    node distance=2.5cm and 0.1cm,
    box/.style={
        rectangle,
        draw=black,
        minimum height=1.5cm,
        text width=6cm,
        align=left,
        inner sep=10pt
    },
    arrow/.style={
        ->,
        thick,
        shorten >=2pt,
        shorten <=2pt,
        >=stealth
    }
]
% Nodes
\node (top) [box, text width=8cm, align=center] {
\textbf{Section \ref{s:system}:} Master model to be approximated \eqref{fulleq} \\[1ex]
$\bullet$ Higher order non-instantaneous collisions\\[0.5ex]
$\bullet$ System of \underline{infinitely many equations} \\[0.5ex]
};

\node (sec3) [box, xshift=-4.5cm, yshift=-5cm, text width=5.5cm] {
    \textbf{Section \ref{s:analysis}:} System for non-instantaneous collisions \eqref{maineps} \\[1ex] 
    $\bullet$ up to $\mathcal{O}(\varepsilon)$ accurate approximation of \eqref{fulleq}\\[0.5ex]
    $\bullet$ system of \underline{two equations}\\[0.5ex]
    $\bullet$ non-instantaneous binary collisions and instantaneous ternary collisions \\[0.5ex]
    $\bullet$ rigorous analysis\\[0.5ex]
};

\node (sec4) [box, xshift=4.5cm, yshift=-5cm, text width=5.5cm] {
    \textbf{Section \ref{s:scalar}:} Scalar model for non-instantaneous collisions \eqref{f1as}\\[1ex]
    $\bullet$ up to $\mathcal{O}(\varepsilon)$ accurate approximation of \eqref{fulleq} \\[0.5ex]
    $\bullet$ \underline{one equation} \\[0.5ex]
    $\bullet$ non-instantaneous binary collisions and instantaneous ternary collisions \\[0.5ex]
    $\bullet$ delay terms
};

% Centered bottom node
\node (standard) [box, text width=5cm, align=center] at ($(sec3.south)!0.5!(sec4.south) - (0,6cm)$) {
     \textbf{Section \ref{s:analysis}:} Standard kinetic model \eqref{limit}\\[1ex]
     $\bullet$ \underline{one equation}\\[0.5ex]
    $\bullet$ instantaneous binary collisions\\[0.5ex]
};

% Connections
\draw [arrow, dotted] (top) -- (sec3);
\draw [arrow, dotted] (top) -- (sec4);

% Bottom Connections - shifted slightly to the corners to open up the middle space
\draw [arrow] (sec3.south east) -- (standard.north);
\draw [arrow, dotted] (sec4.south west) -- (standard.north);

% Labels
% Top label
\node [align=center, text width=2.5cm, font=\small] at ($(top.south) - (0,1.4cm)$) {
    regime for\\
    short collision\\
    duration $\varepsilon \ll 1$
};

% Bottom label - moved HIGHER (-1.2cm instead of -2.2cm)
\node [align=center, font=\small, text width=3cm] at ($(sec3.south)!0.5!(sec4.south) - (0,1cm)$) {
    Instantaneous\\
    limit\\
    $\eps\to 0$
};

\end{tikzpicture}
\caption{
    Overview of the models presented in this article, a dotted arrow indicates the formal passage from one to the other, while a solid arrow indicates a rigorous limit
    }
\end{figure}
\end{center}

\section{System of non-instantaneously interacting particles}\label{s:system}

We shall present a model, which can be interpreted in terms of different applications of alignment, e.g., myxobacteria \cite{balagam} or liquid crystals \cite{oswald}. Indeed, the assumption of instantaneous collisions is not expected to be generally valid in these applications. This issue is discussed, for instance, in \cite{carrillo2025kinetic} (Remark 3.8), where the authors explicitly refer to the paradox described in \cite{palffy2019paradox}. More generally, the latter shows that even for smooth, strictly convex rigid bodies, the classical instantaneous impact framework may lead to unphysical predictions. This suggests that the limitation is not specific to liquid crystals or myxobacteria, but already arises at the level of strictly convex rigid-body interactions. In particular, the post-collisional behavior cannot always be determined solely through an instantaneous collision law, and a time-resolved interaction model may therefore provide a more appropriate description of the dynamics.

Here we will 
describe the dynamics in terms of opinion formation \cite{toscani2006kinetic}, where the opinion of individuals is represented by the one-dimensional variable $v \in \R$. 
Individuals may participate in discussion groups of arbitrary size, where the effect of the discussion is that all participants gradually approach the average opinion of the group.
Two groups may combine to make a bigger group (this effect includes the possibility of individuals joining a group). On the other hand, a group may split into two smaller groups.
Combination and splitting (i.e., coagulation and fragmentation \cite{ball1990discrete}) are governed by Poisson processes with parameters depending on the sizes of the involved groups.

A group of size $k\in\mathbb{N}$ is characterized by the $k$-tuple $(v_1,\ldots,v_k)\in\R^k$ of opinions of its participants. The distribution of groups of size $k$ at time $t\in\R$
will be described by the density $f_k(v_1,\ldots,v_k,t)\ge 0$. The assumption of indistinguishability of the individuals has the consequence that $f_k$ is invariant under permutations
of $(v_1,\ldots,v_k)$. The family $\{f_1,f_2,\ldots\}$ satisfies the system
\begin{equation}\label{fulleq}
    \begin{split}   
    \partial_t\,f_k\ + \nabla_{(v_1,...,v_k)}\cdot (U_k\,f_k) = 
    & \frac{1}{2}\sum_{j=1}^{k-1} \lambda_{j,k-j} f_j\odot f_{k-j} - \sum_{j = 1}^\infty \lambda_{j,k} f_k \int_{\R^j}f_j^*\, d(v_1^*,\dots,v_j^*) \\
    &+ \sum_{j =1}^\infty \mu_{k, j}\int_{\R^j}f_{k+j} \, d(v_{k+1}, \dots,v_{k+j})- \frac{1}{2}\sum_{j= 1}^{k-1}\mu_{j,k-j} f_k\,,
    \end{split}
\end{equation}
$k\ge 1$, where $f_k^*$ denotes evaluation at $(v_1^*,\dots,v_k^*)$, and where the symmetric tensor product in the first term on the right hand side is defined by
$$
   (f_j\odot f_{k-j})(v_1,\ldots,v_k) := {k\choose j}^{-1} \sum_{c\in C(k,j)} f_j(v_c)f_{k-j}(v_{c'}) \,.
$$
Here $C(k,j)$ denotes the set of $j$-combinations of $\{1,\ldots,k\}$. For each \(c=\{i_1,\ldots,i_j\}\in C(k,j)\), we denote \(c':=\{1,\ldots,k\}\setminus c\), and define \(v_c=(v_{i_1},\ldots,v_{i_j})\). The factors $1/2$ in \eqref{fulleq} correct the fact that in the following sums
every term appears twice. The parameters of the above mentioned Poisson processes are
\begin{itemize}
    \item $\lambda_{i,j}=\lambda_{j,i}\ge 0$ for the rate of coagulation between groups of sizes  $i$ and $j$ (beginning of collision or discussion processes)\,,
    \item $\mu_{i,j}=\mu_{j,i}\ge 0$ for the rate of fragmentation of a group of size $i+j$ into groups of sizes $i$ and $j$ (end of collision or discussion processes).
\end{itemize}
The 'acceleration' fields $U_k\in\R^k$, $k\ge 1$, describe the interaction process within a group of $k$ individuals. As indicated above, an alignment process (or trend to the 
average opinion) is assumed:
\begin{equation}\label{Uk}
    (U_k)_i := \frac{1}{k} \sum_{j=1}^k v_j -v_i \,, \quad i \in \{1,\dots,k\}\,, \qquad k>1\,.
\end{equation}
We also define $U_1=0$, meaning that individuals not in a discussion group do not change their opinions. The absence of a factor in front of $v_i$ indicates that in a 
nondimensionalization the relaxation time of the trend towards the average opinion (assumed independent of the size of the discussion group) has been taken as reference time.

\subsection*{Formal properties -- moments}\label{ss:fromalpropfullequ}

The suitability of the model requires certain formal properties. For example, nonnegativity of the distribution functions $f_1,f_2,\ldots$ is formally preserved, since all terms with
a minus sign on the right hand side of \eqref{fulleq} have a factor $f_k$. 

Since individuals only change their opinion, their total number should be preserved by the dynamics. Denoting the total number of groups of size $k$ by $M_k$ 
and the total number of individuals (or the total mass) by $M$, we have
\begin{equation*}
   M_k = \int_{\R^k}f_k \,d(v_1,...v_k)\,,\quad k\ge 1\,,\qquad  M:=\sum_{k=1}^\infty k \,M_k\,.
\end{equation*}
It turns out that the family $\{M_1,M_2,\ldots\}$ solves a closed infinite system of ODEs:
\begin{equation}\label{Mk-system}
\begin{split}
    \dot{M}_k =& \frac{1}{2} \sum_{j= 1}^{k-1} \lambda_{j,k-j} M_j\,M_{k-j} 
               -  \sum_{j=1}^{\infty} \lambda_{j,k} M_j M_k + \sum_{j= 1}^{\infty} \mu_{k,j} M_{k+j} - \frac{1}{2}\sum_{j= 1}^{k-1} \mu_{j,k-j} M_k \,,\quad k\ge 1\,.
\end{split}
\end{equation}
For the rate of change of the total mass we get (with $k=i+j$ in the first and the last term, with the symmetry of the rate constants, and with a symmetrization)
\begin{equation}\label{totalmass}
\begin{split}
    \dot M =  \sum_{k=1}^\infty k \dot{M}_k =& \frac{1}{2} \sum_{i=1}^\infty \sum _{j= 1}^\infty (i+j)\lambda_{j,i} M_j \, M_i - \sum_{k=1}^\infty  \sum _{j= 1}^\infty k\lambda_{j,k} M_k\, M_j \\
     &+\sum_{k= 1}^\infty \sum_{j= 1}^\infty k\mu_{k,j} \,M_{k+j} - \frac{1}{2}\sum_{i= 1}^\infty  \sum_{j= 1}^\infty (i+j)\mu_{j,i} \,M_{i+j} = 0 \,,
\end{split}
\end{equation}
as expected. This book-keeping result is independent from the choice of the rate constants $\lambda_{j,k}$, $\mu_{j,k}$, and of the interaction fields $U_k$. 

System \eqref{Mk-system} is actually the standard discrete coagulation-fragmentation model \cite{ball1990discrete}. Since coagulation and fragmentation can be seen as a
chemical reaction and its reverse, one might hope for the existence of an equilibrium state, where they are balanced. This question cannot be answered in general, but let us
assume that a unique equilibrium $\{M_k^\infty\}_{k\ge 1}$ exists, which has the correct total mass,
$$
    \sum_{k=1}^\infty kM_k^\infty = M \,,
$$
positive partial masses $M_k^\infty>0$, $k\ge 1$, and also satisfies the \emph{detailed balance} condition
$$
     \lambda_{j,k} M_j^\infty M_k^\infty = \mu_{j,k} M_{k+j}^\infty \,,\qquad j,k\ge 1 \,.
$$
This is only possible under various conditions on the rate constants such as, for example, $\lambda_{j,k}>0$ iff $\mu_{j,k}>0$. For uniqueness, a large enough number of pairs $(\lambda_{j,k}, \mu_{j,k})$ needs not to vanish.

By classical results for \emph{mass-action kinetics} \cite{Horn1972massaction}, the \emph{relative entropy}
$$
    H\left(\{M_k\}_{k\ge 1}| \{M_k^\infty\}_{k\ge 1}\right) := \sum_{k=1}^\infty \left( M_k \log \frac{M_k}{M_k^\infty} - M_k + M_k^\infty\right)
$$
is nonincreasing in time. This can be seen by first rewriting the right hand side of \eqref{Mk-system} in terms of $u_k := M_k/M_k^\infty$, $k\ge 1$, and using detailed balance:
$$
   \dot M_k = \frac{1}{2} \sum_{j=1}^{k-1} \mu_{j,k-j} M_k^\infty (u_j u_{k-j} - u_k) + \sum_{j=1}^\infty \mu_{j,k}M_{j+k}^\infty (u_{j+k} - u_j u_k) \,.
$$
After summation against $\log u_k$ and symmetrization of the second term we obtain
$$
  \frac{d}{dt}  H\left(\{M_k\}_{k\ge 1}| \{M_k^\infty\}_{k\ge 1}\right) = - \frac{1}{2} \sum_{j,k=1}^\infty \mu_{j,k} M_{k+j}^\infty (u_{k+j} - u_j u_k) \log\frac{u_{k+j}}{u_j u_k} \le 0 \,.
$$
It is easily seen that under the constraint that the total mass of $\{M_k\}_{k\ge 1}$ is $M$, the right hand side (the \emph{entropy dissipation})
only vanishes for $M_k=M_k^\infty$, $k\ge 1$. This raises the expectation that $\{M_k(t)\}_{k\ge 1}$ converges to $\{M_k^\infty\}_{k\ge 1}$ as $t\to\infty$. It is a classical result
due to Aizenman and Bak \cite{AizenmanBak} that this is true for the case that $\lambda_{j,k}, \mu_{j,k}$ are independent from $(j,k)$ (when $\{M_k^\infty\}_{k\ge 1}$ can be
computed explicitly).

Since the 
interaction in each discussion group preserves the average opinion, we expect the same for the whole ensemble. The average opinion is given by
\begin{equation*}
   v_\infty := \frac{I}{M} \,,\quad\text{with } I=\sum_{k = 1}^\infty k \, I_k\,, \quad  I_k = \frac{1}{k} \int_{\R^k}\sum_{j=1}^k v_j \,f_k \,d(v_1,...,v_k) = \int_{\R^k}v_1 \,f_k \,d(v_1,...,v_k) \,.
\end{equation*}
The family $\{I_1,I_2,\ldots\}$ of first order moments again solves a closed ODE system (assuming to have solved \eqref{Mk-system}): 
\begin{equation}\label{momentum}
\begin{split}
    \dot{I}_k  &= \sum_{j= 1}^{k-1}\lambda_{j,k-j}\, \frac{j}{k}M_{k-j}I_j - \sum_{j=1}^\infty \lambda_{j,k}\, M_j\,I_k  
    + \sum_{j=1}^\infty \mu_{k,j} I_{k+j} - \frac{1}{2}\sum_{j=1}^{k-1} \mu_{j, k-j}I_k \\ &=: \varphi_k(\{I_j\}_{j\ge 1}, \{M_j\}_{j\ge 1})\,.
\end{split}   
\end{equation}
The first term on the right hand side needs some explanation: When the symmetric tensor product $f_j\odot f_{k-j}$ in \eqref{fulleq} is multiplied by $v_1$ and then integrated,
those terms, where $v_1$ appears in the argument of $f_j$, produce the contribution $I_j M_{k-j}$. The other terms produce $M_j I_{k-j}$. The first case occurs ${k-1\choose j-1}$
times, i.e. with probability ${k-1\choose j-1}/{k\choose j}= \frac{j}{k}$. Therefore we obtain 2 terms, which turn out to be the same after the coordinate change $j\mapsto k-j$, namely 
the first term on the right hand side of \eqref{momentum}.

The derivation of \eqref{momentum} also uses the fact that the interaction within the groups does not contribute:
\begin{equation*}
\begin{split}
    \int_{\R^k}&v_1\nabla_{(v_1,...,v_k)}\cdot (U_k\,f_k) d(v_1,...,v_k) = - \int_{\R^k}(U_k)_1 f_k \,d(v_1,...,v_k)\\
    =& - \frac{1}{k}\sum_{j=1}^k \int_{\R^k}v_j f_k\, d(v_1,...,v_k) +  \int_{\R^k}v_1 f_k \, d(v_1,...,v_k) = 0  \,.\\
\end{split}
\end{equation*}
Similarly to \eqref{totalmass} we obtain $\dot I = 0$, showing that the average opinion $v_\infty$ is constant in time.

We shall give a heuristic argument that all the average opinions $\overline v_k := I_k/M_k$, $k\ge 1$, within groups converge to $v_\infty$ as $t\to\infty$, if the 
group sizes $\{M_k\}_{k\ge 1}$ converge to a detailed-balance equilibrium $\{M_k^\infty\}_{k\ge 1}$. We start by writing the right hand side of \eqref{momentum} in terms of the $\overline v_k$ 
and then approximate it for large $t$, replacing $M_k$ by $M_k^\infty$:
\begin{eqnarray*}
   \dot{I}_k  &=& \sum_{j= 1}^{k-1}\lambda_{j,k-j}\frac{j}{k} M_{k-j}M_j\overline v_j - \sum_{j=1}^\infty \lambda_{j,k}\, M_j\,M_k \overline v_k  
    + \sum_{j=1}^\infty \mu_{k,j} M_{k+j}\overline v_{k+j} - \frac{1}{2}\sum_{j=1}^{k-1} \mu_{j, k-j}M_k \overline v_k \\
    &\approx& \frac{1}{2}\sum_{j=1}^{k-1} \mu_{j, k-j}M_k^\infty \left(\frac{j}{k}\overline v_j + \frac{k-j}{k}\overline v_{k-j} - \overline v_k\right) + \sum_{j=1}^\infty \mu_{k,j} M_{k+j}^\infty(\overline v_{k+j} - \overline v_k)
\end{eqnarray*}
The distance of the average opinions to $v_\infty$ can be measured by the quadratic relative entropy
$$
   \frac{1}{2}\sum_{k=1}^\infty k M_k (\overline v_k - v_\infty)^2 \approx \frac{1}{2}\sum_{k=1}^\infty k \frac{I_k^2}{M_k^\infty} - v_\infty I + \frac{1}{2}v_\infty^2 M \,,
$$
and therefore
$$
   \frac{d}{dt}  \frac{1}{2}\sum_{k=1}^\infty k M_k (\overline v_k - v_\infty)^2 \approx \sum_{k=1}^\infty k\overline v_k \dot I_k 
   \approx - \sum_{i,j=1}^\infty \mu_{i,j} M_{i+j}^\infty i (\overline v_i - \overline v_{i+j})^2 \,.
$$
This suggests that all $\overline v_k(t)$ tend to the same value as $t\to\infty$, which has to be $v_\infty$ by the conservation of the average opinion.
As a consequence of the discussion processes,
it is plausible to expect that not only the average opinions of all discussion groups but also the opinion of each individual approaches $v_\infty$. 
This can be checked by introducing the variance
\begin{equation*}
    V = \sum_{k=1}^\infty kV_k \,,
\end{equation*}
with 
\begin{equation*} 
   V_k = \int_{\R^k}(v_1 - v_\infty)^2 f_k \, d(v_1,...,v_k)  = E_k  - 2v_\infty I_k + v_\infty^2 M_k\,,\qquad
    E_k = \int_{\R^k}v_1^2 f_k \, d(v_1,...,v_k) \,.
\end{equation*}
For $k>1$ the time derivative of $E_k$ contains a contribution from the discussion process:
\begin{equation*}
\begin{split}
    &\int_{\R^k}v_1^2\,\nabla_{(v_1,...,v_k)}\cdot (U_k\,f_k) d(v_1,...,v_k) = - 2\int_{\R^k} v_1(U_k)_1 f_k \,d(v_1,...,v_k)\\
    &= - \frac{2}{k^2}\sum_{i=1}^k\sum_{j=1}^k \int_{\R^k}v_i(v_j-v_i) f_k\, d(v_1,...,v_k) =  \frac{1}{k^2} \sum_{i=1}^k\sum_{j=1}^k \int_{\R^k}(v_i-v_j)^2 f_k \, d(v_1,...,v_k) \\
    &= \frac{k-1}{k} \int_{\R^k}(v_1-v_2)^2 f_k \, d(v_1,...,v_k) =:  \frac{k-1}{k} \widetilde V_k \,.
\end{split}
\end{equation*}
The contributions from coagulation and fragmentation are as in \eqref{momentum}. Therefore
$$
    \dot E_k = \varphi_k(\{E_j\}_{j\ge 1}, \{M_j\}_{j\ge 1}) - \frac{k-1}{k} \widetilde V_k \,,\qquad k\ge 1 \,,
$$
with the definition $\widetilde V_1 := 0$. 

As we have seen for the zeroth and first order moments, it turns out that the moments of any order solve a closed ODE system, recursively depending on the lower order moments.
All second order moments can be represented by $E_k$ and $\widetilde V_k$. The time derivative of $\widetilde V_k$ is given by
\begin{equation*}
\begin{split}
    \dot{\widetilde V}_k =& \frac{1}{2} \sum_{j=1}^{k-1} \lambda_{j,k-j} \left( \left( 1 - \frac{2j(k-j)}{k(k-1)}\right)\widetilde V_j M_{k-j} 
              + \frac{2j(k-j)}{k(k-1)} 2\left(E_j M_{k-j}- I_j I_{k-j}\right)\right) \\
    &  - \sum_{j=1}^\infty \lambda_{j,k} \widetilde V_k M_j + \sum_{j=1}^\infty \mu_{k,j}\widetilde V_{k+j} - \frac{1}{2} \sum_{j=1}^{k-1} \mu_{j,k-j} \widetilde V_k - 2\widetilde V_k\,,
         \qquad k>1\,,
\end{split}
\end{equation*}
where in the first line the coefficients $1 - \frac{2j(k-j)}{k(k-1)}$ and $\frac{2j(k-j)}{k(k-1)}$ are the probabilities that, after splitting $\{v_1,\ldots,v_k\}$ into groups of sizes
$j$ and $k-j$, both $v_1$ and $v_2$ end up in the same subgroup and, respectively, in different ones. The last term results from the discussion process.

Finally, we compute the time derivative of the variance:
\begin{equation*}
  \dot V = \sum_{k= 1}^\infty \,k \dot{E}_k = - \sum_{k=2}^\infty (k-1)\widetilde V_k  \le 0\,.
\end{equation*}
This allows the formal conclusion that asymptotically agreement is reached within each discussion group, which 
completes a heuristic argument for the conjecture
$$
    f_k(v_1,\ldots,v_k,t) \to M_k^\infty \prod_{j=1}^k \delta(v_j-v_\infty) \qquad\text{as } t\to\infty \,,\qquad k\ge 1\,,
$$
if \eqref{Mk-system} has a detailed-balance equilibrium.

\subsection*{Fast collision regime and formal first order non-instantaneous approximation}\label{ss:maineps}

We are interested in the regime of almost instantaneous collisions, scaling their average duration with $\eps \ll 1$. Starting from system \eqref{fulleq} we aim to provide a model, which is formally accurate up to $\mathcal{O}(\eps)$. Rescaling equation \eqref{fulleq} in such a way that collisions are short also requires a large fragmentation rate. For the discussions to still have a significant effect, we also need strong interaction fields. As a consequence we expect that larger discussion groups become less likely. This motivates 
\begin{equation}\label{scaling}
    \mu_{i,j}\rightarrow \eps^{-1}\mu_{i,j}\,, \qquad U_k \rightarrow \eps^{-1} U_k \,,\qquad f_k \rightarrow \eps^{k-1}f_k \,,
\end{equation}
with $\eps \ll 1$. This changes \eqref{fulleq} into
\begin{equation}\label{eq:fkeps}
    \begin{split}   
    \eps\partial_t f_k + \nabla_{(v_1,...,v_k)}\cdot (U_k\,f_k) = 
    & \frac{1}{2}\sum_{j=1}^{k-1} \lambda_{j,k-j} f_j\odot f_{k-j} - \sum_{j = 1}^\infty \eps^j\lambda_{j,k} f_k \int_{\R^j}f_j^*\, d(v_1^*,\dots,v_j^*) \\
    &+ \sum_{j =1}^\infty \eps^j\mu_{k, j}\int_{\R^j}f_{k+j} \, d(v_{k+1}, \dots,v_{k+j})- \frac{1}{2}\sum_{j= 1}^{k-1}\mu_{j,k-j} f_k\,.
    \end{split}
\end{equation}
All these equations have fast dynamics, except for $k=1$, where the $O(1)$-terms vanish and the equation can be divided by $\eps$:
\be\label{eq:f1eps}
   \partial_t f_1 = \sum_{j =1}^\infty \eps^{j-1}\mu_{1, j}\int_{\R^j}f_{j+1} \, d(v_2, \dots,v_{j+1}) - \sum_{j = 1}^\infty \eps^{j-1}\lambda_{j,1} f_1 \int_{\R^j}f_j^*\, d(v_1^*,\dots,v_j^*)
\ee
The ensemble of free individuals can gain only from fragmentation and lose only by coagulation.

The instantaneous limit $\eps\to 0$ is the same as in \cite{kanzler2022kinetic}, but we also include its discussion here for the sake of self consistency. The limiting equations
for $k=1$ and $k=2$ are a closed system:
\begin{equation}\label{lim-sys}\begin{split}
   &\partial_t \overline f_1 = \mu_{1, 1}\int_{\R}\overline f_2 \, dv_2 - \lambda_{1,1} \overline f_1 \int_{\R}\overline f_1^*\, dv_1^* \,,\\
   & \nabla_{(v_1,v_2)}\cdot (U_2\,\overline f_2) =  \frac{1}{2} \lambda_{1,1} \overline f_1\otimes \overline f_1 -\frac{1}{2}\mu_{1,1}\overline f_2
\end{split}\end{equation}
This system can be reduced to an equation for $\overline f_1$ after solving the second equation for $\overline f_2$ by the method of characteristics:
\be\label{f2-lim}
   \overline f_2 = \frac{\lambda_{1,1}}{2} \int_0^\infty S_{2,0}(\sigma)(\overline f_1\otimes \overline f_1) d\sigma \,,
\ee
with the semigroup
\be\label{S20}
   (S_{2,0}(\sigma)h)(v_1,v_2) = e^{(1-\mu_{1,1}/2)\sigma} h(v_1',v_2')
\ee
generated by $\overline f_2\mapsto -\nabla_{(v_1,v_2)}\cdot (U_2\,\overline f_2) -\frac{1}{2}\mu_{1,1}\overline f_2$, and with the \emph{collision rule}
\be\label{coll-rule}\begin{split}
    &v_1' = \Phi_{2,1}^{-\sigma}(v_1,v_2) := \frac{1+e^\sigma}{2}v_1 + \frac{1-e^\sigma}{2}v_2 \,,\\
    &v_2' = \Phi_{2,2}^{-\sigma}(v_1,v_2) := \frac{1-e^\sigma}{2}v_1 + \frac{1+e^\sigma}{2}v_2  \,,
\end{split}\ee
connecting the opinions $(v_1',v_2')$ at the beginning of a discussion between two individuals to the opinions $(v_1,v_2)$ at the end of a discussion of duration $\sigma$.
This relation can be inverted by changing the sign of $\sigma$: 
$$
   v_1 = \Phi_{2,1}^\sigma(v_1',v_2') \,,\qquad v_2 = \Phi_{2,2}^\sigma(v_1',v_2') \,.
$$
Substitution of $\overline f_2$ in the equation for $\overline f_1$ results in a kinetic model for binary collisions of standard form:
\be\label{limit}
   \partial_t \overline f_1 = \lambda_{1,1} \int_{\R} \int_0^\infty b(\sigma) \bigl( e^\sigma (\overline f_1\otimes \overline f_1)' - \overline f_1\otimes \overline f_1\bigr)d\sigma\,dv_2 \,,
\ee
with the probability density 
$$
  b(\sigma) = \frac{\mu_{1,1}}{2} e^{-\sigma\mu_{1,1}/2}
$$ 
for the collision duration, with the determinant $e^\sigma$ of the Jacobian of the collision rule, and with the prime denoting evaluation at the pre-collisional state $(v_1',v_2')$. 

The solution of \eqref{limit} approximates the solution component $f_1$ of \eqref{eq:fkeps} formally up to a $\mathcal{O}(\eps)$-error.
It is our goal to improve this approximation by one order in $\eps$. Since in \eqref{eq:f1eps} $f_2$ occurs at leading order, we also need to approximate $f_2$
up to $\mathcal{O}(\eps)$. In both equations for $f_1$ and $f_2$, the component $f_3$ appears in $\mathcal{O}(\eps)$-terms. Therefore we need a leading order approximation for $f_3$.
Since in the equation for $f_3$ components $f_k$, $k>3$, do not occur at leading order, the first step in the approximation procedure is to ignore discussions with more than 3
participants, i.e. $\lambda_{j,k}=\mu_{j,k} = 0$, $j+k>3$. On the other hand we assume that discussion groups with 2 and 3 participants are started and finished:
\begin{equation}\label{rate-ass}
  \lambda_{1,1}, \mu_{1,1}, \lambda_{1,2}, \mu_{1,2} >0 \,.
\end{equation}
This leads to the system
\begin{subequations}\label{scaledsystem}
	\begin{equation}\label{eqf1}
        \p_t f_1 = \mu_{1,1}\int_{\R}f_2\,dv_2 + \eps\,\mu_{1,2}\int_{\R^2}f_3 \,d(v_2,v_3) - \lambda_{1,1}f_1\int_{\R}f_1^*\,dv_1^* 
        - \eps\lambda_{1,2}f_1 \int_{\R^2}f_2^*\,d(v_1^*,dv_2^*)\,,
    \end{equation}
    \begin{equation}\label{eqf2}
        \eps\p_t f_2 + \nabla_{(v_1,v_2)} \cdot (U_2\,f_2) = \frac{1}{2}\lambda_{1,1}f_1\otimes f_1 - \eps \lambda_{1,2}f_2\int_{\R}f_1^* \, dv_1^* + \eps \mu_{1,2}\int_{\R}f_3 \,dv_3
        - \frac{1}{2}\mu_{1,1}f_2 \,,
    \end{equation}
    \begin{equation}\label{eqf3}
        \eps\p_t f_3 + \nabla_{(v_1,v_2,v_3)}\cdot(U_3\,f_3) = \lambda_{1,2} f_1\odot f_2  - \mu_{1,2}f_3 \,,
    \end{equation}
\end{subequations}
Since we only need a leading order approximation of $f_3$, the final approximation step is to consider the quasi-stationary version 
\begin{equation}\label{eqf3qs}
      \nabla_{(v_1,v_2,v_3)}\cdot(U_3\,f_3) = \lambda_{1,2} f_1\odot f_2  - \mu_{1,2}f_3 \,,
\end{equation}
of \eqref{eqf3} and to eliminate $f_3$:
$$
   f_3 = \lambda_{1,2} \int_0^\infty S_3(\sigma)(f_1\odot f_2)d\sigma \,,
$$
where 
$$
   (S_3(\sigma)h)(v_1,v_2,v_3) = e^{(2-\mu_{1,2})\sigma} h(v_1',v_2',v_3')
$$
is the semigroup generated by $f_3 \mapsto -\nabla_{(v_1,v_2,v_3)}\cdot(U_3\,f_3) - \mu_{1,2}f_3$ with the \emph{three-particle collision rule}
 \begin{equation}\label{collrules3}
	\begin{split}
    		v_1' = \Phi_{3,1}^{-\sigma}(v_1,v_2,v_3) := \frac{1+2e^\sigma}{3}v_1 + \frac{1-e^\sigma}{3}v_2 + \frac{1-e^\sigma}{3}v_3 \,, \\
  		v_2' = \Phi_{3,2}^{-\sigma}(v_1,v_2,v_3) := \frac{1-e^\sigma}{3}v_1 + \frac{1+2e^\sigma}{3}v_2 + \frac{1-e^\sigma}{3}v_3\,, \\
  		v_3' = \Phi_{3,3}^{-\sigma}(v_1,v_2,v_3) := \frac{1-e^\sigma}{3}v_1 + \frac{1-e^\sigma}{3}v_2 + \frac{1+2e^\sigma}{3}v_3\,.
    	\end{split}
\end{equation}
Finally we obtain a model, which is formally accurate up to $\mathcal{O}(\eps)$ for both $f_1$ and $f_2$:
\begin{subequations}\label{maineps}
	\begin{equation}\label{mainf1}
	\begin{split}
		\p_t f_1 =& \mu_{1,1} \int_{\R} f_2 \, dv_2 + \eps \mu_{1,2} \lambda_{1,2} \int_{\R^2} \int_0^{\infty}S_3(\sigma) \left( f_1\odot f_2\right) d\sigma \, d(v_2,v_3) \\
		 &- \lambda_{1,1}f_1 \int_{\R}f_1^* \, dv_1^* - \eps \lambda_{1,2} f_1 \int_{\R^2} f_2^* \, d(v_1^*,v_2^*)\,, 
	\end{split}
	\end{equation}
	\begin{equation}\label{mainf2}	
	\begin{split}
		\eps \p_t f_2 + \nabla_{(v_1,v_2)} \cdot \left(U_2 f_2\right)=& \frac{1}{2}\lambda_{1,1} f_1\otimes f_1 -\eps \lambda_{1,2}f_2 \int_{\R} f_1^* \, dv_1^* \\
		& +\eps \mu_{1,2} \lambda_{1,2} \int_{\R}\int_0^{\infty} S_3(\sigma) \left(f_1\odot f_2\right)  d\sigma \, dv_3 
		 - \frac{1}{2}\mu_{1,1} f_2\,,
	\end{split}
	\end{equation}
which will be considered below subject to initial conditions 
\begin{equation}\label{mainICs}
	f_1(v_1,0)=f_1^I(v_1)\,, \quad \quad f_2(v_1, v_2,0)=f_2^I(v_1,v_2)\,,
\end{equation}
\end{subequations}
with the initial data satisfying
\begin{equation}\label{ICs}
    f_1^I,\, f_2^I \ge 0\, , \hspace{0.3cm}\int_{\R}\big(1+v_1^2)f_{1}^I\,dv_1<\infty\, , \hspace{0.3cm}\int_{\R}\big(1+v_1^2)f_{2}^I\,dv_1\,dv_2<\infty\,, \hspace{0.3cm}f_2^I(v_1, v_2) = f_2^I(v_2,v_1).
\end{equation}

\section{Well-posedness and instantaneous limit for the first order accurate model}\label{s:analysis}

This Section \ref{s:analysis} is dedicated to an investigation of model \eqref{maineps}. It contains results on the long-term dynamics of moments
(Subsection \ref{ss:moments}), on existence and uniqueness of solutions (Subsection \ref{ss:existence}), and on the rigorous instantaneous limit (Subsection \ref{ss:limit}).

\subsection{Dynamics of the moments}\label{ss:moments}
 
 We start by deriving formal properties similarly to Section \ref{ss:fromalpropfullequ}. We expect that \eqref{maineps} conserves the total mass
\begin{align}\label{M}
	M:= M_1 + 2\eps M_2 \,, \quad \quad \text{where} \quad M_1:=\int_{\R}f_1 \, dv_1\,, \quad M_2:=\int_{\R^2} f_2 \, d(v_1,v_2)\,.
\end{align}
Note that there is no contribution from $f_3$, since discussions with 3 participants have vanishing duration by \eqref{eqf3qs}, which also implies
$\mu_{1,2}M_3 = \lambda_{1,2}M_1 M_2$. Therefore the partial masses $M_1, M_2$, satisfy the ODE system
\begin{equation}\label{massODE}
	\begin{split}
		\dot{M}_1 =& \mu_{1,1}M_2 - \lambda_{1,1}M_1^2\,,\\
		2\eps \dot{M}_2 = & \lambda_{1,1}M_1^2 - \mu_{1,1}M_2\,,
	\end{split}
\end{equation}
immediately implying the mass conservation
\begin{align*}
	M= M_1(0) + 2\eps M_2(0)\,,
\end{align*}
which can be used to establish convergence $\left(M_1(t), M_2(t)\right) \to \left(M_1^\infty, M_2^\infty\right)$ as $t \to \infty$, with
\begin{align}\label{masses_eq}
	M_1^\infty = \frac{2M}{1+\sqrt{1+ 8\lambda_{1,1}\eps M/\mu_{1,1}}} \,,\qquad M_2^\infty = \frac{\lambda_{1,1}}{\mu_{11}}(M_1^\infty)^2 \,.
\end{align}
The same is expected for the total first moment of the system 
\begin{align}\label{I}
	I:= I_1 + 2\eps I_2 \,, \quad \quad \text{where} \quad I_1:=\int_{\R}v_1\, f_1 \,dv_1\,, \quad I_2:=\int_{\R^2} v_1\,f_2 \, d(v_1,v_2)\,,
\end{align}
where the partial first moments again satisfy a closed ODE system:
\begin{equation}\label{firstmomentODE}
    \dot{I}_1 = - 2\varepsilon\dot{I}_2 =  \mu_{1,1}\,I_2 -\lambda_{1,1}\,M_1\,I_1 + \frac{2}{3} \varepsilon\lambda_{1,2} \left(M_1I_2 -M_2I_1\right)\,,
\end{equation}
where we have again used \eqref{eqf3qs} to get $\mu_{1,2} I_3= \lambda_{1,2}(I_1M_2 + 2 I_2 M_1)/3$. as expected, the first moment is conserved:
$I = I_1(0) + 2 \varepsilon I_2(0)$. Using this for reduction to a scalar equation as well as the convergence of the partial masses, it is obvious that 
\beqar
   I_1(t) &\to& I_1^\infty := \frac{3\mu_{1,1} + 2\eps \lambda_{1,2}M_1^\infty}{3\mu_{11} + 2\eps (3\lambda_{1,1}M_1^\infty + \lambda_{1,2}M)} I = M_1^\infty v_\infty\,,\\
   I_2(t) &\to& I_2^\infty := \frac{3\lambda_{1,1}M_1^\infty + 2\eps \lambda_{1,2}M_2^\infty}{3\mu_{1,1} + 2\eps\lambda_{1,2}M_1^\infty} I_1^\infty = M_2^\infty v_\infty\,,
\eeqar
as $t\to\infty$, with $v_\infty := I/M$, which is independent of time. 

Analogously to the previous section, we define the total variance $V:= V_1 + 2\eps V_2$, where
\begin{align*}
	V_k = \int_{\R^k} (v_1-v_\infty)^2 f_k \,d(v_1,\ldots,v_k) \,,\qquad k=1,2,3\,.
\end{align*}
The partial variances satisfy the ODEs
\begin{eqnarray}
	\dot{V}_1 &=&  \mu_{1,1}V_2 + \varepsilon \mu_{1,2} V_3 - \lambda_{1,1}\,M_1\,V_1 -\varepsilon\lambda_{1,2}M_2\,V_1\,, \nonumber\\
   2\varepsilon \dot{V}_2 &=& \lambda_{1,1} M_1V_1 + 2\eps\mu_{1,2} V_3 - 2\varepsilon\lambda_{1,2}M_1V_2 - \mu_{1,1}V_2 - 4\widetilde V_2\,,\label{V2}
\end{eqnarray}
with
\begin{align*}
	\widetilde V_k = \int_{\R^k} (v_1-v_2)^2 f_k \,d(v_1,\ldots,v_k) \,,\qquad k=2,3\,.
\end{align*}
An equation for $V_3$ is obtained from \eqref{eqf3qs}, after using the computation
\begin{align*}
  &\int_{\R^3}(v_1-v_\infty)^2\nabla_{(v_1,v_2,v_3)}\cdot(U_3 f_3)d(v_1,v_2,v_3) = -\frac{2}{3}\int_{\R^3} v_1 (v_2-v_1 + v_3-v_1)f_3 d(v_1,v_2,v_3) \\
  & = -\frac{4}{3}\int_{\R^3} v_1 (v_2-v_1)f_3 d(v_1,v_2,v_3) = \frac{2}{3} \widetilde V_3 \,,
\end{align*}
where the last equality follows from symmetrization. This implies
\be\label{V3}
   \mu_{1,2} V_3 = \frac{\lambda_{1,2}}{3}(V_1 M_2 + 2 V_2 M_1) - \frac{2}{3}\widetilde V_3 \,.
\ee
We conclude that the variance is nonincreasing:
\be\label{V-decay}
    \dot V = -4\widetilde V_2 - 2 \widetilde V_3 \le 0 \,.
\ee
Further information can be derived from an equation for $\widetilde V_3$, also obtained from \eqref{eqf3qs}:
$$
   \mu_{1,2} \widetilde V_3 = \frac{\lambda_{1,2}}{3}(2V_1 M_2 + 2 V_2 M_1 - 4(I_1-v_\infty M_1)(I_2-v_\infty M_2) + M_1 \widetilde V_2) \,,
$$
implying
$$
    \dot V \le -\frac{4\lambda_{1,2}}{3\mu_{1,2}}(V_1 M_2 + V_2 M_1) + \frac{8\lambda_{1,2}}{3\mu_{1,2}} (I_1-v_\infty M_1)(I_2-v_\infty M_2)  \,.
$$
By our previous results, $M_1$ and $M_2$ converge to positive values as $t\to\infty$, and the last term converges to zero. Therefore for $t$ large enough there exists $\gamma>0$, such that
$$
    \dot V \le - \gamma V + \frac{8\lambda_{1,2}}{3\mu_{1,2}} (I_1-v_\infty M_1)(I_2-v_\infty M_2)  \,,
$$
implying $V(t) \to 0$ as $t\to\infty$ and thus, at least formally,
$$
     f_1(v_1,t) \to M_1^\infty \delta(v_1-v_\infty) \,,\qquad f_2(v_1,v_2,t) \to M_2^\infty \delta(v_1-v_\infty)\delta(v_2-v_\infty) \,,\qquad\text{as } t\to\infty \,.
$$

The above computations can be summarized in the following lemma.
\begin{lemma}[Moment dynamics]\label{lem:moment-dynamics}
System \eqref{maineps} conserves the total mass
$
M=M_1+2\eps M_2
$
and the total first moment
$
I=I_1+2\eps I_2.
$
Moreover, the partial masses and partial first moments satisfy
\[
(M_1(t),M_2(t))\to (M_1^\infty,M_2^\infty),
\qquad
(I_1(t),I_2(t))\to (M_1^\infty v_\infty, M_2^\infty v_\infty),
\qquad \text{as } t\to\infty,
\]
where
\[
v_\infty:=\frac{I}{M},
\qquad
M_1^\infty=\frac{2M}{1+\sqrt{1+8\lambda_{1,1}\eps M/\mu_{1,1}}},
\qquad
M_2^\infty=\frac{\lambda_{1,1}}{\mu_{1,1}}(M_1^\infty)^2\,.
\]
In addition, the total variance
$
V(t)=V_1(t)+2\eps V_2(t)
$
converges to zero as $t \to \infty$.
\end{lemma}

\subsection{Existence and uniqueness}\label{ss:existence}

We start by stating the \emph{mild formulation} of the initial value problem \eqref{maineps}, which can be obtained by integration of the system with respect to time. 
\begin{definition}\label{def:mild_solution}
The pair
\[
(f_1,f_2)\in C\bigl([0,\infty);L^1_+(\mathbb R)\times L^1_+(\mathbb R^2)\bigr)
\]
is called a \underline{global mild solution} of \eqref{maineps}, if it satisfies, for $t\ge 0$,\begin{subequations}\label{mainepsmild}
	\begin{equation}\label{eqf1mild}
	\begin{split}
		f_1(t) =& \, S_{1,\eps}(0,t) f_1^I +  \mu_{1,1} \int_0^t S_{1,\eps}(s,t) \int_{\R} f_2(s) \, dv_2 \, ds \\
		&+ \eps \mu_{1,2} \lambda_{1,2}  \int_0^t S_{1,\eps}(s,t) \int_{\R^2} \int_0^{\infty}S_3(\sigma) \left(f_1(s)\odot f_2(s)\right) d\sigma \, d(v_2,v_3) \, ds\,,
	\end{split}
	\end{equation}
	\begin{equation}\label{eqf2mild}
	\begin{split}
	f_2(t)=&\,S_{2,\eps}\left(0,\frac{t}{\eps}\right)f_2^I + \frac{\lambda_{1,1}}{2\eps}\int_0^t S_{2,\eps}\left(\frac{s}{\eps},\frac{t}{\eps}\right)(f_1(s)\otimes f_1(s)) ds \\
	&+ \mu_{1,2}\lambda_{1,2} \int_0^t  S_{2,\eps}\left(\frac{s}{\eps},\frac{t}{\eps}\right) \int_{\R}\int_0^{\infty} S_3(\sigma) \left( f_1(s)\odot f_2(s)\right) d\sigma \, dv_3  \, ds \,,
	\end{split}
	\end{equation}
\end{subequations}
\end{definition}
In the above definition, \(S_{1,\varepsilon}\) denotes the one-particle semigroup
\begin{align*}
	S_{1,\varepsilon}(s,t) : = \exp\left(-\lambda_{1,1} \int_s^t M_1(r) \, dr -\eps \lambda_{2,1}\int_s^t M_2(r) \, dr\right)\,,
\end{align*}
where 
$$
	M_1:=\int_{\R}f_1 \, dv_1\,, \quad M_2:=\int_{\R^2} f_2 \, d(v_1,v_2)\,,
$$
while \(S_{2,\varepsilon}\) denotes the two-particle semigroup, written in terms of the fast variables \(\sigma=s/\eps\) and \(\tau=t/\eps\),
\begin{equation}\label{d:S2}
	\begin{split}
	 &\big( S_{2,\eps}(\sigma,\tau) h\big)(v_1,v_2) \\ &:=
		\exp\left((1-\mu_{1,1}/2)(\tau-\sigma)-\lambda_{1,2} \int_{\eps\sigma}^{\eps\tau} M_1(r)\, dr \right) h\left(\Phi_{2,1}^{\sigma-\tau}(v_1,v_2),\Phi_{2,2}^{\sigma-\tau}(v_1,v_2)\right)\,,
	\end{split}
\end{equation}
with the coordinate transformations as in \eqref{coll-rule}. We also recall the three-particle semigroup
\[
\bigl(S_3(\sigma)g\bigr)(v_1,v_2,v_3)
:=
e^{(2-\mu_{1,2})\sigma}
g\bigl(\Phi^{- \sigma}_{3,1}(v_1,v_2,v_3),\Phi^{- \sigma}_{3,2}(v_1,v_2,v_3),\Phi^{- \sigma}_{3,3}(v_1,v_2,v_3)\bigr),
\]
as introduced in Section \ref{s:system}. 

Note that the notation is consistent with the previous section in the sense that \eqref{S20} is obtained with 
$\tau-\sigma$ replaced by $\sigma$ and with $\eps=0$.
The arguments $v_1,v_2,v_3$ are suppressed in \eqref{mainepsmild} and in the following, whenever their choice is unambiguous.

\begin{theorem}\label{existence}
	Let $f_1^I \in L_+^1(\R)$, $f_2^I \in L_+^1(\R^2)$, and \eqref{rate-ass} hold. Then \eqref{maineps} has a unique global mild solution.
\end{theorem}
\begin{proof}
Following our considerations regarding the moments in Section \ref{ss:moments}, $M_1$ and $M_2$ can be completely characterized by the dynamics of \eqref{massODE} and therefore can be assumed to be given in the definitions of $S_1$ and $S_{2,\eps}$. It will be used in the following that the semigroups $S_1$, $S_{2,\eps}$, and $S_3$ are
$L^1$-contractions, since the factors $e^{\tau-\sigma}$ in $S_{2,\eps}$ and $e^{2\sigma}$ in $S_3$ are the determinants of the Jacobians of the coordinate transformations
$\Phi_2^{\sigma-\tau}$ and, respectively, $\Phi_3^{-\sigma}$.

Local existence and uniqueness will be proven by Picard iteration. The right hand-side of \eqref{mainepsmild} defines the fixed-point operator 
$\mathcal{F}(f_1,f_2) = (\mathcal{F}_1(f_1,f_2),\mathcal{F}_2(f_1,f_2))$, which obviously preserves positivity 
and, by the contraction property of the semigroups, maps $C([0,T];\,L^1_+(\R) \times L_+^1(\R^2))$ into itself for every $T>0$. More precisely, with the natural norm $\|\cdot\|_{n,T}$ 
on $C([0,T];\,L_+^1(\R^n))$, 
\beqar
  \|\mathcal{F}_1(f_1,f_2)\|_{1,T} &\le& M_1^I + T(\mu_{1,2}\|f_2\|_{2,T} + \eps\lambda_{1,2}\|f_1\|_{1,T} \|f_2\|_{2,T}) \,,\\
  \|\mathcal{F}_2(f_1,f_2)\|_{2,T} &\le& M_2^I + T\left(\frac{\mu_{1,1}}{2\eps}\|f_1\|_{1,T}^2 + \lambda_{1,2}\|f_1\|_{1,T} \|f_2\|_{2,T}\right) \,,\\
\eeqar
with $M_n^I := \int_{\R^n} f_n^I \,d(v_1,\ldots,v_n)$. Here we have used that actually 
\be\label{S3-est}
    \|S_3(\sigma)\|_{L^1(\R^3)\to L^1(\R^3)} \le e^{-\mu_{1,2}\sigma} \,.
\ee
The above estimate implies immediately that, for $T$ small enough, $\mathcal{F}$ maps the set
$$
   \mathcal{S} := \{(f_1,f_2)\in C([0,T];\,L^1_+(\R) \times L_+^1(\R^2)):\, \|f_n\|_{n,T} \le 2M_n^I ,\, n=1,2\}
$$
into itself.

In order to show the contraction property of $\mathcal{F}$ we consider $(f_1, f_2), (\tilde{f}_1, \tilde{f}_2) \in \mathcal{S}$ and show Lipschitz continuity of the second and third terms on the right hand sides of \eqref{eqf1mild} and \eqref{eqf2mild}. The first term
is linear:
    \begin{align*}
         \mu_{1,1} \int_{\mathbb{R}}\bigg| \int_0^t S_1(s,t) \int_{\R} (f_2(s) - \tilde{f}_2(s)) dv_2 \, ds \bigg|\,dv_1 \leq T \mu_{1,1} \|f_2- \tilde{f}_2\|_{2,T}\,,
    \end{align*}
    where we used that $S_1(s,t)\leq 1$. For the second term we again use \eqref{S3-est}:
    \begin{equation*}
        \begin{split}
            &\eps \mu_{1,2} \lambda_{1,2}  \int_{\mathbb{R}}\bigg|\int_0^t S_1(s,t) \int_{\R^2} \int_0^{\infty} S_3(\sigma) \big(f_1(s)\odot f_2(s) - \tilde{f}_1(s)\odot\tilde{f}_2(s)\big) 
            d\sigma \, d(v_2,v_3) ds  \bigg| dv_1 \\
            &\leq \eps \lambda_{1,2} \int_0^t \int_{\mathbb{R}^3} \left|f_1(s)\odot f_2(s)-\tilde{f}_1(s)\odot\tilde{f}_2(s)\right|  d(v_1,v_2,v_3) ds \\
            &\leq \eps \lambda_{1,2} \int_0^t \int_{\mathbb{R}^3} \left(\left|f_1(s)-\tilde{f}_1(s)\right|\odot f_2(s)- \left|f_2(s)-\tilde{f}_2(s)\right|\odot\tilde{f}_1(s)\right) d(v_1,v_2,v_3)  ds \\
            &\leq T\eps \lambda_{1,2} \left( 2M_2^I \|f_1- \tilde{f}_1\|_{1,T} + 2 M_1^I \|f_2- \tilde{f}_2\|_{2,T)} \right) \,.
        \end{split}
    \end{equation*}
 Similar estimates can be carried for the right hand side of \eqref{eqf2mild}. Indeed, we have
    \begin{equation*}
    \begin{split}
    & \frac{\lambda_{1,1}}{2\eps}\int_{\mathbb{R}^2}\left|\int_0^t  S_2\left(\frac{s}{\eps},\frac{t}{\eps}\right)\left(f_1(s)\otimes f_1(s) - \tilde{f}_1(s)\otimes\tilde{f}_1(s)\right)ds \right| 
    d(v_1,v_2) \\
   & \leq T \frac{\lambda_{1,1}}{\varepsilon} 2M_1^I \|f_1- \tilde{f}_1\|_{1,T} \,,
    \end{split}
    \end{equation*}
    and
  \begin{equation*}
    \begin{split}
        & \mu_{2,1}\lambda_{1,2}\int_{\R^2} \left| \int_0^t S_2\left(\frac{s}{\eps},\frac{t}{\eps}\right) \int_{\R}\int_0^{\infty} S_3(\sigma) \left(f_1(s)\odot f_2(s)- \tilde{f}_1(s)\odot\tilde{f}_2(s)\right) d\sigma \, dv_3 \, ds\right| d(v_1,v_2) \\
   &\leq T\lambda_{1,2}  \left( 2M_2^I \|f_1- \tilde{f}_1\|_{1,T} + 2 M_1^I \|f_2- \tilde{f}_2\|_{2,T)} \right) \,.
    \end{split}  
  \end{equation*}
These four estimates show that $\mathcal{F}$ is a contraction for $T$ small enough, implying local existence and uniqueness. 

Our procedure is consistent in the sense that, if $(M_1,M_2)$ is a solution of \eqref{massODE} with $M_n(0)=M_n^I$, n=1,2, and \eqref{mainepsmild} is solved
with this $(M_1,M_2)$ given in the definition of the semigroups, then $\int_{\R^n} f_n d(v_1,\ldots,v_n)= M_n$, $n=1,2$. Proving this statement is a standard procedure, observing that any mild solution of \eqref{maineps} is a weak solution, whence a weak formulation of the moment problem can be derived by an approximation procedure with appropriately chosen families of test functions. Finally, the global existence of $(M_1,M_2)$ implies global existence for \eqref{maineps}.
\end{proof}

\subsection{Instantaneous limit}\label{ss:limit}

The formal instantaneous limit $\eps\to 0$ can be carried out in the same way as in Section \ref{s:system}. As an alternative one may start from \eqref{eqf2mild}, carry out
the coordinate change $s = t-\eps \sigma$, and note that 
$$
    S_{2,\eps}\left(0,\frac{t}{\eps}\right) \to 0 \,,\qquad S_{2,\eps}\left(\frac{t}{\eps}-\sigma,\frac{t}{\eps}\right) \to S_{2,0}(\sigma) \,,\qquad\text{as } \eps\to 0 \,.
$$
Therefore \eqref{f2-lim} is the formal limit of \eqref{eqf2mild}, and the formal limit of $f_1$ satisfies \eqref{limit}. Note that for \eqref{limit} to make sense,
$\lambda_{1,1}, \mu_{1,1}>0$ is needed, which will be assumed from now on. The aim of this section is to make the limit rigorous. 

\begin{theorem}\label{InstantaneousLimit}
Let \eqref{rate-ass} hold and let the initial data satisfy \eqref{ICs} as well as 
$$
   \int_{\R} f_1^I \log(f_1^I)dv_1 + \int_{\R^2} f_2^I \log(f_2^I)d(v_1,v_2) <\infty \,.
$$
Then the mild solution $(f_1, f_2)$ of \eqref{maineps} satisfies
\begin{align*}
        &\lim_{\varepsilon \to 0} f_1(\cdot,t) = \overline f_1(\cdot,t) \hspace{0.3cm} \text{ weakly in } L^1(\R), \text{locally uniformly in } t \in [0,\infty), \\
    &\lim_{\varepsilon \to 0} f_2 = \overline f_2 \quad \text{tightly in }  \R^2\times(0,T), \text{ for any } 0<T<\infty \,,
\end{align*}
where $(\overline f_1,\overline f_2)$ is a weak solution of \eqref{lim-sys} satisfying $\overline f_1(t=0)=f_1^I$.
\end{theorem}

The proof of Theorem \ref{InstantaneousLimit} is based on compactness arguments and will be given after two preliminary steps. First, uniform bounds (as $\eps\to 0$) on the moments will be obtained, which gives compactness in the space of measures. An improvement for $f_1$ to $L^1$-compactness can be achieved in the second step by establishing a uniform bound on the logarithmic entropy. \\

\paragraph{\bf Uniform moment bounds:} The first goal is to obtain uniform-in-$\eps$ bounds for moments of the solution of \eqref{maineps}. It is standard to show (see the last paragraph of the proof of Theorem \ref{existence}) that the moments satisfy the ODEs investigated in Section \ref{ss:moments}.

Starting with the system \eqref{massODE}, the mass conservation property immediately implies $M_1(t)\le M$ with the total mass of the initial data
$$
  M := \int_{\R} f_1^I dv_1 + \eps \int_{\R^2} f_2^I d(v_1,v_2) \,.
$$
Using this in the second equation of \eqref{massODE}, the comparison principle for ODEs immediately implies 
$$
  M_2(t) \le \max\left\{ M_2^I, \frac{\lambda_{1,1}}{\mu_{1,1}}M^2 \right\} \,,
$$
whence $M_1$ and $M_2$ are uniformly bounded as $t\to\infty$ and $\eps\to 0$.

The same procedure is used for the partial variances: A uniform bound for  $V_1$ follows immediately, since the total variance is nonincreasing \eqref{V-decay}. For the second partial variance we use \eqref{V2}
and \eqref{V3} to get
$$
   2\eps \dot V_2 \le \left( \lambda_{1,1}M_1 + \frac{2}{3}\eps \lambda_{1,2}M_2\right)V_1 - \left( \mu_{1,1} + \frac{2}{3}\eps\lambda_{1,2}M_1\right)V_2 \,,
$$
implying uniform boundedness of $V_2$, again by the comparison principle. Collecting these results we have:

\begin{lemma}\label{moments}
    Let \eqref{rate-ass} hold and let the initial data satisfy \eqref{ICs}. Then the solution of \eqref{mainepsmild} satisfies
    \begin{equation*}
        \int_{\R^n} (1+v_1^2)f_n \,d(v_1,\ldots,v_n) < \infty \,,\qquad n=1,2,
    \end{equation*}
uniformly as $\eps\to 0$ and $t\to\infty$.
\end{lemma}

\paragraph{\bf Logarithmic entropy:} 
The leading order coagulation-fragmentation reactions in \eqref{maineps} have the equilibrium $(1,\lambda_{1,1}/\mu_{1,1})$.
Motivated by the theory of chemical reaction networks \cite{Horn1972massaction}, we investigate the corresponding logarithmic relative entropy
\begin{align}
    \mathcal{H}[f_1,f_2] := \int_{\R} f_1 \left(\log(f_1)-1\right)\, dv_1 + \eps \int_{\R^2} f_2 \left(\log\left(\frac{\mu_{1,1}f_2}{\lambda_{1,1}}\right)-1\right) d(v_1,v_2)\,.
\end{align}
Along solutions of \eqref{maineps} we obtain
\begin{align*}
    \frac{d}{dt}\mathcal{H}[f_1,f_2] = & -\frac{1}{2} \int_{\R^2} \left( \lambda_{1,1}f_1\otimes f_1 - \mu_{1,1}f_2\right)\log\left(\frac{\lambda_{1,1}f_1\otimes f_1}{\mu_{1,1}f_2}\right) 
    d(v_1,v_2) \\
    & + \eps \int_{\R^3} \left( \mu_{1,2}f_3 - \lambda_{1,2}f_1\otimes f_2\right)\log\left( \frac{\mu_{1,1}f_1\otimes f_2}{\lambda_{1,1}}\right)d(v_1,v_2,v_3)\\
    & - \int_{\R^2} \nabla_{(v_1,v_2)}\cdot(U_2 f_2) \log\left( \frac{\mu_{1,1}f_2}{\lambda_{1,1}}\right)d(v_1,v_2)\,.
\end{align*}
The first term on the right hand side is the non-positive contribution from the leading order reactions, as expected. The last term, using two integrations by parts as well as
$\nabla_{(v_1,v_2)}\cdot U_2=-1$, can be computed as
$$
  - \int_{\R^2} \nabla_{(v_1,v_2)}\cdot(U_2 f_2) \log\left( \frac{\mu_{1,1}f_2}{\lambda_{1,1}}\right)d(v_1,v_2) = \int_{\R^2} U_2\cdot \nabla_{(v_1,v_2)} f_2\,d(v_1,v_2) = M_2 \,,
$$
which is positive, but uniformly bounded. Finally, with the second term we produce a nonpositive contribution by subtracting and adding
\begin{align*}
   & \eps \int_{\R^3} \left( \mu_{1,2}f_3 - \lambda_{1,2}f_1\otimes f_2\right)\log\left( \frac{\mu_{1,1}\mu_{1,2}f_3}{\lambda_{1,1}\lambda_{1,2}}\right)d(v_1,v_2,v_3) \\
   &= \eps \int_{\R^3} \left( \mu_{1,2}f_3 - \lambda_{1,2}f_1\odot f_2\right)\log\left( \frac{\mu_{1,1}\mu_{1,2}f_3}{\lambda_{1,1}\lambda_{1,2}}\right)d(v_1,v_2,v_3) \\
   & = -\eps \int_{\R^3} \nabla_{(v_1,v_2,v_3)}\cdot (U_3 f_3)\log\left( \frac{\mu_{1,1}\mu_{1,2}f_3}{\lambda_{1,1}\lambda_{1,2}}\right)d(v_1,v_2,v_3) \\
   & = \eps \int_{\R^3} U_3\cdot \nabla_{(v_1,v_2,v_3)} f_3 \,d(v_1,v_2,v_3) = \eps M_3 = \eps \frac{\lambda_{1,2}}{\mu_{1,2}} M_1 M_2\,,
\end{align*}
another positive but uniformly bounded contribution. For the second equality we have used \eqref{eqf3qs}. Combining these results, we have
\begin{align*}
    \frac{d}{dt}&\mathcal{H}[f_1,f_2] \le M_2 + \eps \frac{\lambda_{1,2}}{\mu_{1,2}} M_1 M_2 \,.
\end{align*}
This formal computation can be justified by first considering a regularized version of $\mathcal{H}[f_1,f_2]$ and subsequent passage to the limit. We also observe that the assumptions of Theorem \ref{InstantaneousLimit} imply $\mathcal{H}[f_1^I,f_2^I] < \infty$.

\begin{lemma}\label{lem:entropy} 
Let the assumptions of Theorem \ref{InstantaneousLimit} hold and let $0<T<\infty$. Then $f_1 \log(f_1)$ is bounded in $L^\infty((0,T), L^1(\R))$ uniformly as $\eps\to 0$.
\end{lemma}

\paragraph{\bf Passing to the limit:} The instantaneous limit will be based on the uniform bounds in Lemmas \ref{moments} and \ref{lem:entropy}.

\begin{proof}[Proof of Theorem \ref{InstantaneousLimit}]
Let $0<T<\infty$. Then Lemma \ref{moments} implies that $\{f_1\}_\eps $ and $\{f_2\}_\eps $ are tight sets of measures on $\R\times(0,T)$ and, respectively, $\R^2\times(0,T)$. 
Due to the Prokhorov theorem \cite{prokhorov1956convergence} this is equivalent to weak sequential compactness of $\{f_1\}_\eps $ and $\{f_2\}_\eps $ in the spaces of positive measures $\mathcal{M}^+(\R\times(0,T))$ and, respectively, $\mathcal{M}^+(\R^2\times(0,T))$. The uniform tightness of $f_1$, together with the uniform entropy bound in Lemma \ref{lem:entropy}, allow to apply the Dunford-Pettis criterion and to deduce weak sequential compactness of $\{f_1\}_\eps$ in $L^1(\R \times (0,T))$. Additionally, from \eqref{mainf1} one can obtain the following $L^1$-bound on the time derivative of $f_1$:
\begin{equation}\label{boundderivativef1}
    \left\|\partial_t f_1 \right\|_{L^1(\R)} \leq \max\{\mu_{1,1} M_2,\lambda_{1,1} M_1^2\} + \eps\lambda_{1,2}M_1 M_2 \,,
\end{equation}
implying (by Lemma \ref{moments}) uniform Lipschitz continuity of the map $t \mapsto f_1(\cdot,t)$ with respect to the $L^1(\R)$-topology. Hence, we can deduce the existence of an accumulation point $f_1^0 \in C([0,\infty); L^1(\R))$ of the family $\{f_1\}_\eps$, such that for a sequence $\varepsilon_n \to 0$, the sequence $\{f_1\}_{\eps_n}$
converges to $f_1^0$ with respect to $L^1(\R)$, locally uniformly in $t \in [0,\infty)$. 

For the fast variable $f_2$ there is no uniform bound of the time derivative as for $f_1$. Therefore we only obtain tight convergence (up to subsequences).

Passing to the limit in the weak formulation of \eqref{maineps} is straightforward, since the $\mathcal{O}(\eps)$-terms tend to zero by uniform boundedness, and the only leading order nonlinearity is $f_1\otimes f_1$, where we use that weak convergence of two measures implies weak convergence of the product measure to the product measure of the limits (\cite{billingsley2013convergence}, Theorem 2.8 (ii)). Finally, the restriction to subsequences is not necessary, since uniqueness for the initial value problem for \eqref{limit} can be shown analogously to the proof of 
Theorem \ref{existence}.
\end{proof}

\section{A first order accurate, non-instantaneous, scalar model}\label{s:scalar}

In Section \ref{ss:maineps} we derived system \eqref{maineps}, describing a kinetic model including first order non-instantaneous correction terms for the interaction processes of the particles. This system promises to be a good candidate to model particle dynamics with close to instantaneous interactions, as the analysis in Section \ref{s:analysis} shows. An obvious further question, matter of discussion in this section, is whether one can find a scalar equation, which gives a well-posed first order non-instantaneous approximation of the dynamics. 

The basic idea will be to eliminate $f_2$ from \eqref{maineps}. This is made feasible by the observation that in the first term on the second line of \eqref{mainf2}, $f_2$ can be 
replaced by its formal limit as $\eps\to 0$, since this will only introduce an $O(\eps^2)$-error. The same argument could be used for the last term in the first line, but this would 
obstruct non-negativity of the approximation. Therefore we start from the mild formulation \eqref{eqf2mild}, where we introduce the coordinate transformation $s=t-\eps\sigma$:
\begin{equation}\label{eqf2mild1}\begin{split}
	f_2(t)=&\,S_{2,\eps}\left(0,\frac{t}{\eps}\right)f_2^I + \frac{\lambda_{1,1}}{2}\int_0^{t/\eps} S_{2,\eps}\left(\frac{t}{\eps}-\sigma,\frac{t}{\eps}\right)
	(f_1(t-\eps\sigma)\otimes f_1(t-\eps\sigma)) d\sigma \\
	&+ \eps\mu_{1,2}\lambda_{1,2} \int_0^{t/\eps}  S_{2,\eps}\left(\frac{t}{\eps}-\sigma,\frac{t}{\eps}\right) \int_{\R}\int_0^{\infty} S_3(\rho) 
	\left( f_1(t-\eps\sigma)\odot f_2(t-\eps\sigma)\right) d\rho \, dv_3  \, d\sigma \,.
\end{split}\end{equation}
The right hand side will be approximated by dropping the first term, since it is exponentially small away from $t=0$, and by passing to the limit in the coefficient of $\eps$ in the
second line:
\begin{equation}\label{f2as}\begin{split}
	f_{2,\text{as}}(t)=&\, f_{2,1}(t) + \eps f_{2,2}(t) \\
	:=&\, \frac{\lambda_{1,1}}{2}\int_0^{t/\eps} \exp\left(-\lambda_{1,2}\int_{t-\eps\sigma}^t M_1(s)ds\right)S_{2,0}(\sigma)
	(f_1(t-\eps\sigma)\otimes f_1(t-\eps\sigma)) d\sigma \\
	&+ \eps\mu_{1,2}\lambda_{1,2} \int_0^\infty  S_{2,0}(\sigma) \int_{\R}\int_0^{\infty} S_3(\rho) 
	\left( f_1(t)\odot f_2^0(t)\right) d\rho \, dv_3  \, d\sigma \,,
\end{split}\end{equation}
with 
$$
    f_2^0(t) = \frac{\lambda_{1,1}}{2}\int_0^\infty S_{2,0}(\sigma)(f_1(t)\otimes f_1(t)) d\sigma \,.
$$
Note that \eqref{f2as} is a formal approximation of $f_2$ with an $O(\eps^2)$-error, given explicitly in terms of $f_1$. Now these approximations are used in \eqref{mainf1}
($f_{2,\text{as}}$ in the leading order terms and $f_2^0$ in the $\mathcal{O}(\eps)$-terms):
\begin{equation}\label{f1as}
  \p_t f_1 = Q_2(f_1,f_1) + \eps Q_3(f_1,f_1,f_1) \,,
\end{equation}
with
\begin{align*}
&Q_2(f_1,f_1) = \lambda_{1,1} \int_{\R} \Biggl[ \frac{\mu_{11}}{2}\int_0^{t/\eps} \exp\left(-\lambda_{1,2}\int_{t-\eps\sigma}^t M_1(s)ds\right) \\
& \hskip 5cmS_{2,0}(\sigma)(f_1(t-\eps\sigma)\otimes f_1(t-\eps\sigma)) d\sigma - f_1\otimes f_1\Biggr] dv_2 
\end{align*}
and
$$
Q_3(f_1,f_1,f_1) =  \lambda_{1,2} \int_{\R^2} \left[ \mu_{1,2} \left( 1 + \mu_{1,1}\int_0^\infty S_{2,0}(\sigma)d\sigma\right)\int_0^\infty S_3(\rho)(f_1\odot f_2^0)d\rho
- f_1\otimes f_2^0\right] d(v_2,v_3) \,.
$$
In the binary collision operator $Q_2$ collisions are non-instantaneous, causing delays in the gain term, as already observed in \cite{kanzler2022kinetic}.
The ternary collision operator $Q_3$ is an instantaneous approximation. Taking into account non-instantaneous effects from the full model would only create 
$O(\eps^2)$-corrections. The gain term of $Q_3$ involves iterated applications of the two-particle and three-particle semigroups, since a ternary collision requires a predecessing
binary collision to happen. This structure is somewhat reminiscent of the Wild sum representation \cite{wild1951boltzmann} of solutions of the Boltzmann equation,
which has been related to iterated higher order collisions, e.g., by Villani \cite{villani2002review}. However, different from that our semigroups also contain an account of the 
dynamics during collisions.

Model \eqref{f1as} preserves non-negativity of $f_1$. However, it does not conserve mass, which is no surprise because $f_1$ does not represent the particles involved in 
non-instantaneous binary collisions. This has also been observed in models for non-instantaneous collisions of quantum particles \cite{lipavsky1999noninstantaneous}, where an 
auxiliary \emph{correlated density} is introduced as a correction. In the following it will be shown that the one-particle marginal of $2 \eps f_{2,1}$ in \eqref{f2as} is a good approximation
for the correlated density. Note that the contributions to $f_{2,\text{as}}$ satisfy the equations
\begin{eqnarray*}
   \eps\p_t f_{2,1} + \nabla_{(v_1,v_2)} (U_2 f_{2,1}) &=& \frac{\lambda_{1,1}}{2} f_1\otimes f_1 - \left(\frac{\mu_{1,1}}{2} + \eps\lambda_{1,2}M_1\right)f_{2,1} \,,\\
   \nabla_{(v_1,v_2)} (U_2 f_{2,2}) &=& \mu_{1,2} \int_{\R} f_3 \,dv_3 - \frac{\mu_{1,1}}{2} f_{2,2} \,,\qquad\text{with} \\
   \nabla_{(v_1,v_2,v_3)} (U_3 f_3) &=& \lambda_{1,2} f_1\odot f_2^0 - \mu_{1,2}f_3 \,,
\end{eqnarray*}
implying
\begin{eqnarray*}
 && \eps \dot M_{2,1} = \frac{\lambda_{1,1}}{2} M_1^2 - \left(\frac{\mu_{1,1}}{2} + \eps\lambda_{1,2}M_1\right) M_{2,1} \,,\\
 && \frac{\mu_{1,1}}{2} M_{2,2} = \mu_{1,2} M_3 = \lambda_{1,2}M_1 M_2^0 \,.
\end{eqnarray*}
These and integration of \eqref{f1as} give
$$
    \frac{d}{dt} \left( M_1 + 2\eps M_{2,1}\right) = 2\eps \lambda_{1,2} M_1 \left(M_2^0 - M_{2,1}\right) = O(\eps^2) \,,
$$
which is the desired result up to an $O(\eps^2)$-error.

\vspace{1cm}

 \subsection*{Acknowledgments:} L.K. received funding by the European Commission under the Horizon2020 research and innovation programme, Marie Sklodowska-Curie grant agreement No 101034255. C.M and C.S. acknowledge support from the Austrian Science Fund, grant numbers W1261 and SFB65.

\subsection*{Data Availibility.}
 No datasets were generated or analysed during the current study.

\subsection*{Declaration of Interests.}
 All authors certify that they have no affiliations with or involvement in any organization or entity with any financial interest or non-financial interest in the subject matter or materials discussed in this manuscript.


\begin{thebibliography}{90}

\bibitem{AizenmanBak} M. Aizenman and T.A. Bak. \emph{Convergence to Equilibrium in a System of Reacting Polymers.} Communications in Mathematical Physics 65 (1979), 
pp. 203--230.

\bibitem{ampatzoglou2021rigorous} I. Ampatzoglou and N. Pavlović. \emph{Rigorous derivation of a ternary Boltzmann equation for a classical system of particles.} Communications in Mathematical Physics 387 (2021), no. 2, pp. 793--863.

\bibitem{aranson2005pattern} I.S. Aranson and L.S. Tsimring. \emph{Pattern formation of microtubules and motors: Inelastic interaction of polar rods.} Physical Review E 71 (2005), no. 5, pp. 050901.

\bibitem{balagam} R. Balagam and O.A. Igoshin. \emph{Mechanism for Collective Cell Alignment in Myxococcus xanthus Bacteria.} PLOS Comp. Biol. 11(8) (2015), e1004474.

\bibitem{ball1990discrete} J.M. Ball and J. Carr. \emph{The discrete coagulation-fragmentation equations: existence, uniqueness, and density conservation.} Journal of Statistical Physics 61 (1990), pp. 203-234.

\bibitem{bertin2009hydrodynamic} E. Bertin, M. Droz, and G. Grégoire. \emph{Hydrodynamic equations for self-propelled particles: microscopic derivation and stability analysis.} Journal of Physics A: Mathematical and Theoretical 42 (2009), no. 44, 445001.

\bibitem{billingsley2013convergence} P. Billingsley. \emph{Convergence of probability measures.} John Wiley \& Sons, 2013. 

\bibitem{bobylev2000some} A.V. Bobylev, J.A. Carrillo, and I.M. Gamba. \emph{On some properties of kinetic and hydrodynamic equations for inelastic interactions.} Journal of Statistical Physics 98 (2000), pp. 743--773.

\bibitem{bolley2012mean} F. Bolley, J.A. Ca\~{n}izo, and J.A. Carrillo. \emph{Mean-field limit for the stochastic Vicsek model.} Applied Mathematics Letters 25 (2012), no. 3, pp. 339--343. 

\bibitem{briant2022cauchy} M. Briant, A. Diez, and S. Merino-Aceituno. \emph{Cauchy theory for general kinetic Vicsek models in collective dynamics and mean-field limit approximations.} SIAM Journal on Mathematical Analysis 54 (2022), no. 1, pp. 1131--1168. 

\bibitem{calvez2014confinement} V. Calvez, G. Raoul, and C. Schmeiser. \emph{Confinement by biased velocity jumps: aggregation of \emph{Escherichia coli}.}  Kinetic and Related Models 8 (2015), no. 5, pp. 651--666.

\bibitem{carlen2015boltzmann} E. Carlen, M.C. Carvalho, P. Degond, and B. Wennberg. \emph{A Boltzmann model for rod alignment and schooling fish.} Nonlinearity 28 (2015), no. 6, pp. 1783.

\bibitem{carrillo2025kinetic} J.A. Carrillo, P.E. Farrell, A. Medaglia, and U. Zerbinati. \emph{A Kinetic Theory Approach to Ordered Fluids.} arXiv preprint (2025), arXiv:2508.10744.
  
\bibitem{carrillo2007contractive} J.A. Carrillo and G. Toscani. \emph{Contractive probability metrics and asymptotic behavior of dissipative kinetic equations.} Rivista di Matematica della Università di Parma, Series 7, Volume 6 (2007), pp. 75--198.


\bibitem{cercignani2013mathematical} C. Cercignani, R. Illner, and M. Pulvirenti. \emph{The mathematical theory of dilute gases.} Springer Science \& Business Media, 2013. 

\bibitem{driver1973equation} R.D. Driver, D.W. Sasser, and M.L. Slater. \emph{The equation \( x'(t) = ax(t) + bx(t - \tau) \) with “small” delay.} The American Mathematical Monthly 80 (1973), no. 9, pp. 990--995.

\bibitem{during2015opinion} B. Düring and M.-T. Wolfram. \emph{Opinion dynamics: inhomogeneous Boltzmann-type equations modelling opinion leadership and political segregation.} Proceedings of the Royal Society A 471 (2015), no. 2182, pp. 20150345.

\bibitem{harvey2011study} C.W. Harvey, F. Morcos, C.R. Sweet, D. Kaiser, S. Chatterjee, X. Liu, D.Z. Chen, and M. Alber. \emph{Study of elastic collisions of \emph{Myxococcus xanthus} in swarms.} Physical Biology 8 (2011), no. 2, pp. 026016.

\bibitem{hittmeir2020kinetic} S. Hittmeir, L. Kanzler, A. Manhart, and C. Schmeiser. \emph{Kinetic modelling of colonies of myxobacteria.} Kinetic and Related Models 14 (2021), pp. 1-24.

\bibitem{Horn1972massaction} F. Horn, R. Jackson, General mass action kinetics, Archive Rat. Mech. Anal. 47 (1972), pp.  81-116.

\bibitem{kanzler2022kinetic} L. Kanzler, C. Schmeiser, and V. Tora. \emph{Two kinetic models for non-instantaneous binary alignment collisions.} Kinetic and Related Models 17 (2024), no. 5, pp. 697-712.

\bibitem{khalil2002nonlinear} H.K. Khalil. \emph{Nonlinear Systems}. Third edition, Prentice Hall, Upper Saddle River, NJ, 2002.

\bibitem{kurzweil2006small} J. Kurzweil. \emph{Small delays don't matter.} In: \emph{Proceedings of the Symposium on Differential Equations and Dynamical Systems}, pp. 47--49, Springer, 2006.

\bibitem{lipavsky1999noninstantaneous} P. Lipavský, V. Špička, and K. Morawetz. \emph{Noninstantaneous collisions and two concepts of quasiparticles.} Physical Review E 59 (1999), no. 2, pp. R1291.

\bibitem{mischler2006cooling} S. Mischler, C. Mouhot, and M. Rodriguez Ricard. \emph{Cooling process for inelastic Boltzmann equations for hard spheres, Part I: The Cauchy problem.} Journal of Statistical Physics 124 (2006), pp. 655--702.

\bibitem{MURPHY2024109266} P. Murphy, M. Perepelitsa, I. Timofeyev, M. Lieber-Kotz, B. Islas, and O.A. Igoshin. \emph{Breakdown of Boltzmann-type models for the alignment of self-propelled rods.} Mathematical Biosciences 376 (2024), pp. 109266.

\bibitem{oswald} P. Oswald and P. Pieranski. \emph{Nematic and Cholesteric Liquid Crystals: Concepts and Physical Properties Illustrated by Experiments},
CRC Press, 2005.

\bibitem{palffy2019paradox} P. Palffy-Muhoray, E.G. Virga, M. Wilkinson, and X. Zheng. \emph{On a paradox in the impact dynamics of smooth rigid bodies.} Mathematics and Mechanics of Solids 24 (2019), no. 12, 4215--4222.

\bibitem{prokhorov1956convergence} Y.V. Prokhorov. \emph{Convergence of random processes and limit theorems in probability theory.} Theory of Probability \& Its Applications 1 (1956), no. 2, pp. 157--214.

\bibitem{sengers1973three} J.V. Sengers. \emph{The three-particle collision term in the generalized Boltzmann equation.} In: \emph{The Boltzmann Equation: Theory and Applications}, pp. 177--208, Springer, 1973.

\bibitem{toscani2006kinetic} G. Toscani. \emph{Kinetic models of opinion formation.} Communications in Mathematical Sciences 4 (2006), no. 3, pp. 481--496.

\bibitem{verhulst2006methods} F. Verhulst. \emph{Methods and applications of singular perturbations: boundary layers and multiple timescale dynamics.} Springer, 2006.

\bibitem{vicsek1995novel} T. Vicsek, A. Czirók, E. Ben-Jacob, I. Cohen, and O. Shochet. \emph{Novel type of phase transition in a system of self-driven particles.} Physical Review Letters 75 (1995), pp. 1226--1229.

\bibitem{villani2002review} C. Villani. \emph{A review of mathematical topics in collisional kinetic theory.} Handbook of Mathematical Fluid Dynamics, vol. 1, pp. 71--305, 2002.

\bibitem{wild1951boltzmann} E. Wild. \emph{On Boltzmann's equation in the kinetic theory of gases.} Mathematical Proceedings of the Cambridge Philosophical Society 47 (1951), no. 3, pp. 602--609.

\end{thebibliography}
\end{document}